\documentstyle[11pt]{article}

\newenvironment{procl}[1]{\medskip\par\noindent {\bf #1.}  \it}%
{\par\noindent}

\title{Entire functions and compact operators with $S_p$-imaginary
component}

\author{Vladimir Matsaev
\thanks
{Supported in part by the Israel Science Foundation of the Israel Academy
of Sciences and Humanities under Grant No. 93/97-1 and by the United
States -- Israel Binational Science Foundation under Grant No. 96-00030. 
}
\
and Mikhail Sodin$^{\ast}$}

\date{}

\begin{document}

\maketitle
\begin{abstract}
We study operators of the form $X+Y$ where $Y$ has a finite $p$-th
Schatten norm, $1\le p <2$, and $X$ is selfadjoint and of Hilbert-Schmidt
class. Our study is based on new theorems on zero distribution of entire
functions of finite order.
\end{abstract}

\subsection{Introduction}

Let $S_p$ be the von Neumann-Schatten
class of compact operators in a Hilbert space $\cal H$ such that
$$
||A||_p = \left\{ \sum_k s_k^p(A) \right\}^{1/p} <\infty\,,
$$ 
where $\{s_k(A)\}$ is a sequence of singular values of the operator $A$,
that is eigenvalues of a non-negative operator $\sqrt{A^*A}$, enumerated
in non-increasing order: $s_1(A)\ge s_2(A)\ge ...\,$. If $p\ge 1$, this
defines a norm, and for these $p$'s the classes $S_p$ are Banach spaces.    
In particular, $S_1$ is the trace class, and $S_2$ is the Hilbert-Schmidt
class. 

For any operator $A$, we denote by 
$$
G=\frac{1}{2}(A+A^*) \qquad {\rm and} \qquad H=\frac{1}{2i}(A-A^*)
$$
its Hermitian parts. 

If $A\in S_2$, then we define the Carleman determinant of $A$ as a
canonical product of genus one 
$$
C_A(z) = \prod_k E(z\mu_k(A))\,, \qquad E(z)=(1-z)e^z\,,
$$
where $\{\mu_k(A)\}$ is the set of eigenvalues of $A$ counted their
multiplicities.

If $A=G+iH$ is a compact quasinilpotent operator, then, by a theorem of
L.~Sakhnovich \cite{sakh}, $||G||_2=||H||_2$, and $||A||_2\le 
2\,||H||_2$. Soon after Sakhnovich obtained his result, M.~Krein
discovered in \cite{krein} a weak type inequality for compact
quasinilpotent operators with imaginary component of trace class (see
inequality (1.6) below), and the first-named author proved
in \cite{matsaev} that, for $1<p<\infty$, $||A||_p\le K_p\,||H||_p$.
This can be regarded as boundedness of the operator of
triangular integration (which recovers a compact quasinilpotent operator 
by its imaginary part) in the spaces $S_p$, $1<p<\infty$, and its ``weak
boundedness'' in the space $S_1$ (cf. \cite{branges}, \cite{gokr2}, and  
\cite[Chapter~4]{davidson}). There are numerous 
interplays and interrelations between  these results and inequalities of
M.~Riesz and Kolmogorov for conjugate harmonic functions (cf.
\cite{matsaev}, \cite[Chapter~4]{gokr2}, \cite{muhly}, and
\cite{sukochev}).   

Original proofs given in \cite{krein} and \cite{matsaev} were based on the
entire function theory. Later, more transparent geometrical proofs have
been found. Recently, working on \cite{ms3}, we realized that function
theory methods work in a more general situation, when  the absolute value
of the Carleman determinant $C_A$ is not too small on the real axis. (If
the operator $A$ is quasinilpotent, then $C_A\equiv 1$.) This required a
further development of function theory methods which was done in
\cite{ms1}. Operator theory results of the present work may be regarded as
non-commutative counterparts of those of \cite{ms1}. Probably, these
results are new even for finite matrices. Our approach also leads to new
results on entire functions (Theorems~3, 5 and 7 below) which may be
interesting on their own.

\medskip
We thank Iossif Ostrovskii for very helpful correspondence and remarks.
 
\subsection{1. Results and conjectures}

First, introduce some notations. 
For an entire function $f(z)$, we denote 
$$
M(r,f)=\max_{|z|=r} |f(z)|\,,
$$
then 
$$
\sigma (f) = \limsup_{r\to\infty} \frac{\log M(r)}{r}
$$
is the exponential type of $f$. By $n(r,f)$ we denote the counting
function of zeros of $f(z)$: it is the number of zeros of $f(z)$ in the
disk $\{|z|\le r\}$ counting multiplicities.  

An entire function $f(z)$ is said to be of
Cartwright class, if its exponential type is finite ($0\le \sigma (f)
<\infty$), and 
$$
\int_{-\infty}^\infty \frac{\log^+|f(t)|}{t^2+1}\, dt <\infty\,.
$$

Further, we set 
$$
\alpha_p(f) = \frac{1}{\pi} \int_{-\infty}^{\infty} 
\frac{\log^+|f(t)|}{|t|^{p+1}}\,dt\,, \qquad \alpha(f)=\alpha_1(f)\,,
$$
by $\alpha_p(f^{-1})$ and $\alpha (f^{-1})$ we denote similar integrals
with $\log^-|f|$ instead of $\log^+|f|$. We also set  
$$
\gamma_p(f)=\sum_k \left| {\rm Im}\frac{1}{z_k} \right|^p\,, \qquad
\gamma(f) = \gamma_1(f)\,,
$$
where the sum is taken over the zero set of $f$ (counting multiplicities).
We allow to these quantities to be equal to $+\infty$. 

By $K$ we denote various positive numerical constants, by $K_p$ we denote
positive values which may depend on $p$ only. 

\begin{procl}{Theorem 1} 
If $A\in S_2$, then
$$
\sigma(C_G)+\alpha(C_G) = \alpha(C_G^{-1}) \le
||H||_1+\alpha(C_A^{-1})\,.
\eqno (1.1)
$$
\end{procl}

In particular, if the RHS of (1.1) is finite, then the Carleman
determinant $C_G(z)$ is of Cartwright class.  
Using a theory of entire functions of Cartwright class (cf. \cite{levin},
\cite{koosis}), we obtain an additional information about distribution of 
eigenvalues $\{\mu_k(G)\}$ of $G$. For example, 
\begin{itemize}
\item[$\bullet$] for each $s>0$,
$$
{\rm card}\, \{k:\, |\mu_k(G)|\ge s \} \le \frac{K}{s}\left(||H||_1 +
\alpha(C_A^{-1})\right)\,;
\eqno (1.2)
$$
\item[$\bullet$] if the RHS of (1.1) is finite, then there exist the 
limits 
$$
\lim_{k\to\pm\infty} k\mu_k(G) = \frac{2}{\pi} \sigma(C_G)\,.
\eqno (1.3)
$$ 
Here we assume that $\mu_k$ are numerated in the decreasing order, and
the positive (negative) values of indices $k$ correspond to the positive
(negative) values of $\mu_k$.
\end{itemize}

Applying estimate (1.2) and Lemma~5.1 from
\cite[Chapter~III]{gokr2}, we obtain a weak-type inequality for singular
values of $A$:
\begin{eqnarray*}
{\rm card }\, \{k: s_k(A) &\ge& s \} \\
&\le& {\rm card }\, \{k: |\mu_k(G)|\ge s/2 \} + {\rm card }\, \{k:
|\mu_k(H)|\ge s/2 \} \\
&\le& 
\frac{K}{s}\left(||H||_1 + \alpha(C_A^{-1})\right)\,.
\qquad \qquad \qquad \qquad \qquad \quad \ (1.4)
\end{eqnarray*}
Now, we mention two special cases. First, assume that the
operator $A$ is quasinilpotent, that is ${\rm spec }(A)=\{0\}$.  Then
$C_A(z)\equiv 1$, and $\alpha (C_A^{-1})=0$. By a theorem of L.~Sakhnovich
\cite{sakh} (see also \cite[Chapter~I]{gokr2}), a quasinilpotent compact
operator $A$ with an imaginary component of trace class belongs to $S_2$. 
We arrive at

\begin{procl}{Theorem 2} If $A$ is a compact quasinilpotent operator with
imaginary component of trace class, then
$$
\sigma(C_G)+\alpha(C_G) = \alpha(C_G^{-1}) \le
||H||_1\,,
\eqno (1.5)
$$
and
$$
{\rm card }\, \{k: s_k(A) \ge s \} \le \frac{K}{s}||H||_1 \,.
\eqno (1.6)
$$ 
\end{procl}

This result, in essence, coincides with a theorem of M.~Krein
\cite[Chapter~IV, Theorem~8.2]{gokr1} which states that $C_G(z)$ is of
Cartwright class, $\sigma(C_G)\le ||H||_1$, and relation (1.6) holds. Our
estimate (1.5) is more precise, and this is important for some
applications (cf. \cite{ms3}). Our proof of Theorem~1 uses Krein's
arguments. We also mention that for quasinilpotent operators, Krein
computed the type $\sigma(C_G)$ using triangular truncation (cf.
\cite[Chapter~VII]{gokr2}). 

\medskip
An interesting question appears here which we were unable to answer: 
{\it Let $\Pi(z)$ be  an arbitrary
canonical product of genus one with real zeros and of Cartwright class,
whether there exists a compact quasinilpotent operator $A=G+iH$ with a
trace class imaginary component $H$ and such that $C_G=\Pi$?}

\medskip
Now, we consider another special case. In the space $l^2=l^2({\bf N})$ we 
define the operator of multiplication
$$
(A_w\,a)(k)=w(k)a(k)\,, \qquad \lim_{k\to\infty} w(k)=0\,.
$$
Then $A_w$ is a compact operator, and
its spectrum coincides with the set $\{w(k)\}$. The singular values of
$A_w$ are $\{|w(k)|\}$. If $w\in l^2$, then $A_w\in S_2$, and 
$$
C_{A_w}(z)=\prod_k E(zw_k)
$$
is an arbitrary canonical product of genus one.

\begin{procl}{Theorem 3}
Let $\Pi (z)$ be a canonical product of genus one. Then 
$$
\sigma (\Pi) + \alpha (\Pi) = 
2\max\left\{\gamma_+(\Pi),\gamma_-(\Pi)\right\}+\alpha(\Pi^{-1}) \,,
\eqno (1.7)
$$
where $\gamma_{+}(\Pi)$
(correspondingly, $\gamma_-(\Pi)$) is 
the sum $\sum\, |\hbox{Im}\,(z_k^{-1})|$ taken over zeros of $\Pi$
lying in the upper (lower) half-plane.
\end{procl}

Equation (1.7) states that if one of its sides is finite, than the other
is finite as well, the both sides are equal, and in particular $\Pi (z)$
belongs to the Cartwright class.  We also get a sharp estimate
$$
\sigma (\Pi) + \alpha (\Pi) \le  
2\gamma (\Pi)+\alpha(\Pi^{-1}) \,.
\eqno (1.8)
$$ 
Estimate (1.8) yields uniform upper bounds for 
$\log M(r,\Pi)$ and for the counting function of zeros $n(r,\Pi)$ in terms
of $\gamma(\Pi)+\alpha(\Pi^{-1})$ 
(cf \cite[Theorem~2]{ms1}).

\medskip
The next result deals with compact operators with imaginary
component from the class $S_p$, $p>1$. It is motivated by a
result of
the first-named author \cite{matsaev} which states that {\it if an
imaginary
component of a compact quasinilpotent operator $A$ belongs to $S_p$,
$1<p<\infty$, then $A\in S_p$, and
$$
||A||_p \le K_p ||H||_p\,.
\eqno (1.9)
$$
}
Another motivation comes from our recent work \cite{ms1}.

\begin{procl}{Theorem 4}
Let $A\in S_2$. Then, for $1<p<2$,
$$
||A||_p\le K_p \left(||H||_p + \alpha_p^{1/p}(C_A^{-1}) \right)\,.
\eqno (1.10)
$$
\end{procl}

In the case $p=2$, one almost immediately gets 
$$
||A||_2 \le 2||H||_2\,,
$$
provided that $A\in S_2$, and 
$$
\int_0 \frac{\log^-|C_A(t)|}{|t|^3}\, dt < \infty\,.
$$
Indeed, we have
$$
\log|C_A(t)| = -\frac{t^2}{2} \mbox{Re}\,\mbox{tr}A^2 + O(t^3)
=\frac{t^2}{2}(\mbox{tr}H^2-\mbox{tr}G^2) + O(t^3)\,, \qquad t\to 0\,,
$$
and since the integral over a neighbourhood of the origin
converges, 
the coefficient by $t^2$ must be non-negative: 
$$
\mbox{tr}(H^2)\ge \mbox{tr}(G^2)\,,
$$
whence
$$
||A||_2\le ||G||_2+||H||_2 \le 2||H||_2\,.
$$

\medskip
As special cases of Theorem~4, we obtain estimate (1.9) for compact
quasinilpotent operators (in this case, by a duality argument
\cite{gokr2}, \cite{davidson}, the case $p>2$ follows from the case
$1<p<2$); and the following result on entire functions:

\begin{procl}{Theorem 5}
Let $\Pi (z)$ be a canonical product of genus one.  Then, for $1<p<2$,
$$
\int_0^{\infty} \frac{\log M(r,\Pi)}{r^{p+1}}\,dr
\le K_p \left(\gamma_p(\Pi) + \alpha_p(\Pi^{-1}) \right)\,.
\eqno (1.11)
$$
\end{procl}

In reality, we first prove Theorem~5 and then, with its help, we establish
Theorem~4.

\medskip
A precise value of the best possible
constant $K_p$ in the inequality
$$
||G||_p\le K_p ||H||_p\,, \qquad 1<p<\infty\,,
$$
valid for compact quasinilpotent operators is still unknown for the most
of values of $p$. There is a long-standing conjecture that it coincides
with the value of the best possible constant in the inequality of
M.~Riesz found by Pichorides:
$$
K_p=\left\{ \begin{array}{ll}
\tan (\pi/2p) & \mbox{if $1<p\le 2$} \\ \\ 
\cot (\pi/2p) & \mbox{otherwise}\,.
\end{array}
\right.
$$
This is confirmed only for $p=2^n$ and $p=2^n/(2^n-1)$, $n\in{\bf N}$
\cite{krupnik}.
A value of the  best possible constant in a more general inequality (1.10)    
is also unknown. 

\medskip
We cannot extend Theorems~4 to other normed ideals of compact
operators (cf. \cite{gokr1}). For example, a natural conjecture is that,
{\it for $1<p<2$,
$$
\sup_{n\in {\bf N}} \left\{s_n(A)n^{1/p} \right\} \le K_p
\sup_{n\in {\bf N}} \left\{s_n(H)n^{1/p}\right\}
$$
provided that $A\in S_2$ and $|C_A(x)|\ge 1$ for $x\in {\bf R}$.}

\medskip

In Theorem~4, we assumed that $\alpha_p(C_A^{-1})<\infty$, that is the
function $x\mapsto \log|C_A(x)|$ is ``almost non-negative'' on the real
line. In this form, the theorem cannot be extended to $p>2$, as an example
suggested by Nazarov shows (cf. \cite{ms1}). However, there is another
possible way to extend inequality (1.9) to more general classes of compact
operators than quasinilpotent ones. Following Levin \cite{levin1}, we
define a regularized determinant $D_A(z)$ of an arbitrary compact operator
$A$ as an entire function whose zero set coincides with the sequence
$\{\mu_k^{-1}(A)\}$, and such that $D_A(0)=1$. Naturally, such a
choice is
not unique. Now, we require that the function $\log^-|D_A(z)|$ is
relatively
small in  $\bf C$:
$$
\int_0^\infty \frac{m(r,D_A^{-1})}{r^{p+1}}\,dr <\infty\,,
$$
where 
$$
m(r,f)=\frac{1}{2\pi} \int_0^{2\pi} \log^+|f(re^{i\theta})|\, d\theta
$$ 
is the Nevanlinna proximity function.  

\begin{procl}{Theorem 6}
Let $A$ be a compact operator, and let $D_A$ be a
regularized determinant of $A$. Then, for $1<p\le 2$,
$$
||A||_p^p \le K_p\, \left(||H||_p^p +
\int_0^\infty \frac{m(r,D_A^{-1})}{r^{p+1}}\,dr
\right)\,.
$$
\end{procl}

Observe, that we do not require here that $A$ is a Hilbert-Schmidt class
operator.
The proof of this result is similar to that of Theorem~4, but is more
lengthy, and we will not give it in this paper. We are unable to
extend Theorem~6 to the case $p>2$, such an extension seems plausible. 
However, the corresponding result from the entire function theory holds
for the whole range $1<p<\infty$ (cf. \cite[Theorem~5]{ms1}):
\begin{procl}{Theorem 7}
Let $f(z)$ be an entire function such that $f(0)=1$. Then for
$1<p<\infty$,
$$
\int_0^\infty \frac{\log M(r,f)}{r^{p+1}}\,dr
\le K_p \, \left( \gamma_p(f)+
\int_0^\infty \frac{m(r,1/f)}{r^{p+1}}\,dr
\right).
$$
\end{procl}

\medskip The last remark is that it seems to be important to have
``prototheorems'' which would contain the results of this work and our
previous work \cite{ms1} as special cases. A natural context for such
prototheorems is provided by von Neumann operator algebras with finite or
semi-finite traces (cf. \cite{dixmier}). For estimates (1.9) and (1.6)
such generalizations are known (cf. \cite{muhly} and \cite{sukochev}). 
In \cite{brown}, L. Brown introduced subharmonic counterparts of the
regularized determinants for operators which belong to one of
non-commutative $L^p(M,\tau)$ classes, where $M$ is a von Neumann algebra
and $\tau$ is a faithful, normal, semi-finite trace on $M$. Using his
results and the arguments given below, one can prove a rather general
version of Theorem~1. Probably, a similar version of Theorem~4 also holds.   

\subsection{2. Proofs of Theorems~1 and 3}

\par\noindent{\it Proof of Theorem~1:} Let $P_{\pm}$ be the orthogonal
projectors onto subspaces generated by those of eigenfunctions of $H$
which correspond to non-negative (negative) eigenvalues. Let 
$H_{\pm}=\pm HP_{\pm}$. The operators $H_{\pm}$ are non-negative, and
$H=H_+-H_-$. Set $H_1=H_+ + H_-$. Then eigenvalues of $H_1$ are absolute
values of eigenvalues of $H$ and $||H||_1={\rm tr}(H_1)$.

Introduce an auxiliary operator $A_1=G+iH_1$. Since $G-A\in S_1$ and
$A-A_1\in S_1$, we can define the perturbation determinants: 
$$
\Delta_{G/A_1}(z)={\rm det}\left[(I-zG)(I-zA_1)^{-1}\right]
=\frac{C_G(z)}{C_{A_1}(z)}\, e^{iz\,{\rm tr}(H_1)}\,,
\eqno (2.1)
$$
and 
$$
\Delta_{A_1/A}(z)={\rm det}\left[(I-zA_1)(I-zA)^{-1}\right]
=\frac{C_{A_1}(z)}{C_A(z)}\, e^{iz\,{\rm tr}(H-H_1)}
\eqno (2.2)
$$  
(cf. \cite[Chapter IV]{gokr1}).
Combining (2.1) and (2.2), we obtain
$$
C_G(z) = \Delta_{G/A_1}(z)\,\Delta_{A_1/A}(z)\, C_A(z)\,  e^{-iz\,{\rm
tr}(H)}\,.
\eqno (2.3)
$$

A chief fact here is a result which goes back to M.~Liv\^{s}ic
\cite[Chapter~IV, Theorem~5.2]{gokr1}: {\it if $A=G+iH$ is compact and
dissipative (that is, $H\ge 0$), and ${\rm tr}(H)<\infty$, and if 
$B=G+iF$ is compact and $-H\le F\le H$, then $|\Delta_{B/A}(z)|\le 1$, for
$z$ in the upper half-plane.}

Applying twice this theorem, we obtain
$$
|\Delta_{G/A_1}(z)|\le 1\,, \qquad {\rm Im}z \ge 0\,,
\eqno (2.4)
$$
and
$$
|\Delta_{A_1/A}(z)|=|\Delta_{A/A_1}(z)|^{-1}\ge 1\,, \qquad {\rm Im}z \ge
0\,.
\eqno (2.5)
$$
Thus, by virtue of (2.3), 
$$
\log^-|C_G(t)|\le \log\frac{1}{|\Delta_{G/A_1}(t)|} +\log^-|C_A(t)|\,,
$$
and
$$
\frac{1}{\pi} \int_{-\infty}^\infty \frac{\log^-|C_G(t)|}{t^2}\, dt \le
- \frac{1}{\pi} \int_{-\infty}^\infty
\frac{\log|\Delta_{G/A_1}(t)|}{t^2}\, dt
+ \frac{1}{\pi} \int_{-\infty}^\infty \frac{\log^-|C_A(t)|}{t^2}\, dt\,.
\eqno (2.6)
$$

Now, we estimate the first integral in the RHS of (2.6). Since
$\log|\Delta_{G/A_1}(z)|$ is a non-positive harmonic function in the upper
half-plane,  
$$
\log |\Delta_{G/A_1}(iy)| = cy + \frac{y}{\pi}\int_{-\infty}^\infty
\frac{\log|\Delta_{G/A_1}(t)|}{t^2+y^2}\, dt 
$$
with $c\le 0$. Whence
$$
\frac{y}{\pi}\int_{-\infty}^\infty
\frac{\log|\Delta_{G/A_1}(t)|}{t^2+y^2}\, dt 
\ge \lim_{y\downarrow 0}\frac{\log |\Delta_{G/A_1}(iy)|}{y}\,.
$$
The latter limit can be computed using relation (2.1): since $C_{A_1}$ and
$C_G$ are canonical products of genus one,
$$
\log|C_{A_1}(z)| = O(z^2)\,, 
\qquad {\rm and} \qquad
\log|C_{G}(z)| = O(z^2)\,, 
$$
for $z\to 0$, and 
$$
\lim_{y\downarrow 0}\frac{\log |\Delta_{G/A_1}(iy)|}{y} = -{\rm
tr}(H_1)\,.
$$
Returning back to (2.6), we have
$$
\frac{1}{\pi} \int_{-\infty}^\infty \frac{\log^-|C_G(t)|}{t^2}\, dt \le
||H||_1 + \alpha(C_A^{-1})\,.
\eqno (2.7)
$$

Next, we prove that, for each $\epsilon >0$,
$$
\frac{1}{\pi} \int_{-\infty}^\infty \frac{\log|C_G(t)|}{t^2}\,
e^{-\epsilon t^2}\, dt \le 0\,.
\eqno (2.8)
$$
Indeed, denoting by $\{\xi_k\}$ the zero set of $C_G$, and integrating
termwise the formula
$$
\log|C_G(t)| = \sum_k \log |E(t/\xi_k)|\,,
$$
we obtain
\begin{eqnarray*}
\int_{-\infty}^\infty \frac{\log|C_G(t)|}{t^2}\, 
e^{-\epsilon t^2}\, dt 
&=& \sum_k \int_{-\infty}^\infty
\frac{\log|E(t/\xi_k)|}{t^2}\, e^{-\epsilon t^2}\, dt \\
&=& \sum_k \int_0^\infty \frac{\log|1-t^2/\xi_k^2|}{t^2}\,
e^{-\epsilon t^2}\, dt \\
&=& \sum_k \left[ \int_0^{\sqrt{2}|\xi_k|} + \int_{\sqrt{2}|\xi_k|}^\infty
\right]
\frac{\log|1-t^2/\xi_k^2|}{t^2}\, e^{-\epsilon t^2}\, dt \\
&\le& \sum_k e^{-2\epsilon \xi_k^2}
\int_0^\infty \frac{\log|1-t^2/\xi_k^2|}{t^2}\, dt \\ \\
&=& 0\,,
\end{eqnarray*}
since
$$
\int_0^\infty \frac{\log|1-t^2|}{t^2}\, dt =0\,.
$$
This proves (2.8).

Now, we obtain that, for each $\epsilon >0$,
\begin{eqnarray*}
\frac{1}{\pi} \int_{-\infty}^\infty \frac{\log^+|C_G(t)|}{t^2}\,
e^{-\epsilon t^2}\, dt &\stackrel{(2.8)}\le& 
\frac{1}{\pi} \int_{-\infty}^\infty \frac{\log^-|C_G(t)|}{t^2}\,
e^{-\epsilon t^2}\, dt \\
&\stackrel{(2.7)}\le& ||H||_1 + \alpha(C_A^{-1})\,,
\end{eqnarray*}         
and therefore, by the monotone convergence theorem, 
$$
\alpha(C_G) \le ||H||_1 + \alpha(C_A^{-1})\,.
\eqno (2.9)
$$
It remains to improve a little bit estimate (2.9): in its LHS we must
replace $\alpha (C_G)$ by $\sigma(C_G)+\alpha(C_G)$. 

First, observe that since $C_G(z)$ has a minimal type with respect to
order two, we have, for $z$ in the upper half-plane,
$$
\log|C_G(z)|=\frac{y}{\pi}\int_{-\infty}^\infty
\frac{\log|C_G(t)|}{|t-z|^2}\,dt + c_0 + c_1{\rm Im}z + c_2{\rm Im}z^2\,.
\eqno (2.10)
$$
The Poisson integral in the RHS of (2.10) is absolutely convergent due to
estimates (2.7) and (2.9). Putting in (2.10) $z=iy$, $y\downarrow 0$, we
obtain that $c_0=0$. Putting there $z=re^{i\pi/4}$, $r\to\infty$, we
obtain that $c_2=0$. At last, putting again $z=iy$, dividing the both
sides by $y$, and letting $y\to 0$, we obtain 
$$
0=\frac{1}{\pi} \int_{-\infty}^{\infty} \frac{\log|C_G(t)|}{t^2}\,dt +
c_1\,,
$$ 
or
$$
c_1+\alpha(C_G) = \alpha(C_G^{-1})\,.
$$
On the other hand, we see from (2.10) (with $c_0=c_2=0$) that 
$$
c_1=\lim_{y\to\infty} \frac{\log |C_G(iy)|}{y} = \sigma(C_G)\,.
$$
This completes the proof of Theorem~1. $\Box$

\medskip\par\noindent{\it Proof of Theorem~3:} We keep the previous
notations which are applied now to the operator of multiplication $A=A_w$
defined above (before formulation of Theorem~3). Since we already know
that $C_G(z)$ 
is of Cartwright class, it has a bounded type; i.e. is a quotient of
bounded analytic functions in the upper half-plane. Then, 
using relations (2.3)--(2.5), we see that $\Pi(z)=C_A(z)$ 
also has a bounded type in the upper half-plane. The same argument shows,
that $\Pi (z)$ has a bounded type in the lower half-plane. Applying a
theorem of M.~Krein (cf. \cite[Lecture~16]{levin}), we obtain that
$\Pi (z)$ is of Cartwright class. 

We need some classical facts about entire functions of Cartwright class 
which may be found in \cite[Lectures~16 and 17]{levin}.
Let $f(z)$ be an entire function of Cartwright type. 
The quantities
$$
\sigma_{\pm}(f) = \limsup_{y\to +\infty} \frac{\log |f(\pm iy)|}{y}
$$  
are called the types of $f$ in the upper and lower half-planes
correspondingly. 
Then
\begin{itemize}
\item[(i)] $$\sigma (f) = \max \{\sigma_+(f),\sigma_-(f) \};$$
\item[(ii)] 
$$
\lim_{R\to\infty} \frac{1}{\pi R}\,\int_0^{\pi} 
\log|f(Re^{i\theta})|\sin\theta\,d\theta = \frac{\sigma_+(f)}{\pi}\,
\int_0^{\pi} \sin^2\theta\, d\theta = \frac{\sigma_+(f)}{2}\,;
$$
\end{itemize}

Now, we apply the Carleman integral formula for the upper half-plane (cf.
\cite[Lecture~24]{levin}) to the function $\Pi (z)$. Due to condition
$\log|\Pi(z)|=O(|z^2|)$ as $z\to 0$, it has no error term, and we get the
relation
$$
2\,\sum_{0<\theta_k<\pi}\left( \frac{1}{r_k}-\frac{r_k}{R^2}
\right)^+\sin\theta_k
$$
$$
=\frac{1}{\pi} \int_{-\infty}^{\infty} \left(\frac{1}{t^2}-\frac{1}{R^2}
\right)^+
\log|\Pi(t)|\,dt + \frac{2}{\pi R} \int_0^\pi
\log|\Pi(Re^{i\theta})|\sin\theta\,d\theta\,,
\eqno (2.11)
$$
where $\{z_k\}=\{r_ke^{i\theta_k}\}$ are zeros of $\Pi$.
Using the dominate convergence theorem and property (ii), 
we easily make in (2.11) the limit transition for
$R\to\infty$, and obtain the limit relation
$$
2\gamma_+(\Pi) = 
\frac{1}{\pi}\, \int_{-\infty}^{+\infty} \frac{\log|\Pi(t)|}{t^2}\,dt
+ \sigma_+(f)\,,
$$
or
$$
\sigma_+(\Pi) + \alpha(\Pi) = 2\gamma_+(\Pi) +
\alpha(\Pi^{-1})\,. 
\eqno (2.12)
$$

Applying the same argument in the lower half-plane, we obtain the second
relation
$$
\sigma_-(\Pi) + \alpha(\Pi) = 2\gamma_-(\Pi) +
\alpha(\Pi^{-1})\,. 
\eqno (2.13)
$$
Relations (2.12), (2.13) and property (i) of entire functions of
Cartwright class complete the proof of Theorem~3.
$\Box$

\medskip One can prove Theorem~3 directly, without appeal to
Theorem~1, along similar lines. Then, instead of a deep result from
operator theory due to Liv\^sic, we need its ``scalar version'':
a simple inequality which  
concerns complex numbers: {\it if $w_1, w_2$ are complex numbers such that
${\rm Re}\,w_1={\rm Re}\,w_2$ and $|{\rm Im}w_1|\le {\rm Im}\,w_2$,
then, for all $z$ in the upper half-plane, 
$$
\left| \frac{1-zw_1}{1-zw_2} \right| < 1\,.
$$    
}

\medskip
I.~Ostrovskii recently proposed another (a
little bit shorter) proof of Theorem~3 based on a different idea. 

\subsection{3. Proofs of Theorems 4 and 5}

First, we prove Theorem~5 and then Theorem~4.

\medskip\par\noindent{\it Proof of Theorem~5:} Let $\{z_k\}$ be zeros of
the canonical product $\Pi(z)$ of genus one. We prove that
$$
\sum_k \frac{1}{|z_k|^p} \le K_p
\left( \gamma_p(\Pi)+\alpha_p(\Pi^{-1}) \right)\,.
\eqno (3.1)
$$ 
Then a standard application of Borel's estimate of the canonical product
(cf. \cite[Lecture~4]{levin}) yields estimate (1.11):
$$
\log M(r,\Pi) \le Kr\,
\left(\int_0^r\frac{n(t,\Pi)}{t^2}\,dt 
+r\,\int_r^\infty\frac{n(t,\Pi)}{t^3}\,dt \right)\,,
$$
whence
\begin{eqnarray*}
\int_0^{\infty} \frac{\log M(r,\Pi)}{r^{p+1}}\,dr &\le&
K\,\int_0^\infty\, \frac{dr}{r^p} 
\left(\int_0^r\frac{n(t,\Pi)}{t^2}\,dt +
r\,\int_r^\infty\frac{n(t,\Pi)}{t^3}\,dt
\right) \\ \\ 
&\le& K\int_0^\infty \frac{n(t,\Pi)}{t^2}\,dt 
\left(\frac{1}{t}\int_0^t
\frac{dr}{r^{p-1}}+\int_t^\infty\frac{dr}{r^p}\right) \\ \\
&\le& K_p\int_0^\infty \frac{n(t,\Pi)}{t^{p+1}}\,dt \\ \\ 
&\stackrel{(3.1)}\le&
K_p\left(\gamma_p(\Pi)+\alpha_p(\Pi^{-1}) \right)\,.
\end{eqnarray*}

Let $\theta_k=\arg z_k$, $-\pi\le\theta_k\le \pi$. Since, for
$|\varphi|\le \pi/2$ and $1<p<2$, $\cos p\varphi\le K_p  \cos^p\varphi$,
we have
$$
\sum_k \frac{\cos^+p(|\theta_k|-\pi/2)}{|z_k|^p}
\le K_p\, \sum_k \frac{|\sin\theta_k|^p}{|z_k|^p} = K_p\gamma_p(\Pi)\,.
\eqno (3.2)
$$
  
\begin{procl}{Claim}
$$
\sum_k \frac{\cos^-p(|\theta_k|-\pi/2)}{|z_k|^p}
\le \sum_k \frac{\cos^+p(|\theta_k|-\pi/2)}{|z_k|^p}
+K_p\alpha_p(\Pi^{-1})\,,
\eqno (3.3)
$$
with
$$
K_p= \frac{p\sin\frac{\pi p}{2}}{2\pi}\,.
$$
In particular, the LHS of (3.3) is finite, provided that
$\gamma_p(\Pi)+\alpha_p(\Pi^{-1})<\infty $.   
\end{procl}

First, we show that the claim easily yields inequality (3.1), and then we
shall prove the claim.

Fix $\eta$, $0<\eta<\frac{\pi}{2}-\frac{\pi}{2p}$, and split the sum in
the LHS of (3.1) into two parts. Into the first sum (which we denote
$\sum^{'}$) we include those $k$'s that $\eta\le |\theta_k|\le \pi -\eta$,
and in the second sum (denoted by $\sum^{''}$)
we include such $k$'s that either
$|\theta_k|<\eta$, or $|\theta_k \pm \pi|\le \eta$. Then
$$
{\textstyle \sum_k ^{'}}\, \frac{1}{|z_k|^p} \le K_p\,
{\textstyle \sum_k}\,\left|{\rm Im}\frac{1}{z_k}
\right|^p\,,  
$$
and, using (3.2) and (3.3),
$$  
{\textstyle \sum_k^{''}}\, \frac{1}{|z_k|^p} \le K_p\,
{\textstyle \sum_k}\,
\frac{\cos^-p(|\theta_k|-\pi/2)}{|z_k|^p}
\le K_p\left(\gamma_p(\Pi)+\alpha_p(\Pi^{-1})
\right)\,.
$$
proving (3.1).

\medskip\par\noindent{\it Proof of the Claim:} 
For each $\epsilon>0$, we have 
\begin{eqnarray*}
-\alpha_p(\Pi^{-1}) &\le& 
\int_{-\infty}^\infty \frac{\log|\Pi(x)|}{|x|^{p+1}}\, e^{-\epsilon x^2}\,
dx \\ \\ 
&=& 
\sum_k \int_{-\infty}^\infty \frac{\log|E(x/z_k)|}{|x|^{p+1}}\,
e^{-\epsilon x^2}\, dx \\ \\ 
&=&
\sum_k \int_0^\infty \frac{\log|1-x^2/z_k^2)|}{|x|^{p+1}}\,
e^{-\epsilon x^2}\, dx \\ \\ 
&=& 
\sum_k \frac{1}{|z_k|^{p+1}} \int_0^\infty
\frac{\log|1-se^{-2i\theta_k}|}{s^{p/2+1}}\,
e^{-\epsilon s |z_k|^2}\, ds\,. 
\qquad \qquad \quad (3.4)
\end{eqnarray*}
Now, 
\begin{eqnarray*}
\int_0^\infty \frac{\log|1-se^{-2i\theta}|}{s^{p/2+1}}\,
e^{-\epsilon s\rho^2} \,ds 
&=&
\left(\int_0^{2\cos^+2\theta}+\int_{2\cos^+2\theta}^\infty \right)
\frac{\log|1-se^{-2i\theta}|}{s^{p/2+1}}\,
e^{-\epsilon s\rho^2} \,ds \\ \\
&\le& e^{-2\epsilon \rho^2 \cos^+2\theta}\,
\int_0^\infty \frac{\log|1-se^{-2i\theta}|}{s^{p/2+1}}\,ds\,,
\end{eqnarray*}
since the kernel $s\mapsto \log|1-se^{-2i\theta}|$ is negative for 
$0<s<2\cos^+2\theta$ and positive for $2\cos^+2\theta<s<\infty$. Computing
the latter integral:
$$
\int_0^\infty \frac{\log|1-se^{-2i\theta}|}{s^{p/2+1}}\,ds
=\frac{\pi}{\frac{p}{2}\sin\frac{\pi p}{2}}\,
\cos\frac{p}{2}(2|\theta|-\pi)\,, \qquad -\pi\le\theta\le \pi\,,
$$ 
we get
$$
\int_0^\infty \frac{\log|1-se^{-2i\theta}|}{s^{p/2+1}}\,
e^{-\epsilon s\rho^2} \,ds \le  
\frac{\pi}{\frac{p}{2}\sin\frac{\pi p}{2}}\,
e^{-2\epsilon \rho^2 \cos^+2\theta} \cos p(|\theta|-\pi/2)\,.
\eqno (3.5)
$$
Plugging (3.5) into (3.4) with
$\theta=\theta_k$ and $\rho=|z_k|$, we obtain
$$
\sum_k \frac{\cos^-p(|\theta_k|-\pi/2)}{|z_k|^p}\,
e^{-2\epsilon|z_k|^2\cos 2\theta_k} \le 
\sum_k \frac{\cos^+p(|\theta_k|-\pi/2)}{|z_k|^p} \,+\, 
\frac{p\sin\frac{\pi p}{2}}{2\pi}\,
\alpha_p(\Pi^{-1})\,.
$$
By the monotone convergence theorem applied for $\epsilon\to 0$, this
proves the claim and completes
the proof of the theorem. $\Box$


\medskip\par\noindent{\it Proof of Theorem~4} uses the strategy from
\cite{matsaev} (see also \cite[Chapter~IV]{matsaev2}) combined with
Theorem~5. Define the
resolvent of $G$
$$
R_G(z)=(I-zG)^{-1}\,.
$$
Since $H\in S_p$ and $S_p$ is an ideal, $HR_G(z)\in S_p\subset S_2$, and
therefore, for regular values of $z$,  $[HR_G(z)]^2\in S_1$.
Hence we can define the determinant 
$$
D(z)={\rm det}[I+z^2(HR_G(z))^2]\,.
\eqno (3.6)
$$
Since   
$$
I+z^2[HR_G(z)]^2 = (I-zA^*)(I-zG)^{-1}(I-zA)(I-zG)^{-1}\,,
\eqno (3.7)
$$
we obtain the factorization 
$$
D(z) = \frac{C_A(z)C_{A^*}(z)}{C_G^2(z)}
\eqno (3.8)
$$
which plays an important role in our considerations. Formula (3.8) follows
immediately, if $A\in S_1$; in the general case, we first approximate
$A\in S_2$ by finite rank operators $A_n$, verify relation (3.8) for
$A_n$, and then let $n\to\infty$. 

Now, using relation (3.6), we derive a good upper bound for
$|D(re^{i\theta})|$, $\theta\ne 0,\pi$. According to one of the
equivalent definitions of the determinant (cf. \cite[Chapter~IV]{gokr1}),
we have 
$$
D(z) = \prod_k \left(1+z^2\mu_k^2(HR_G)\,\right)\,,
$$  
where $\mu_k(HR_G)$ are eigenvalues of $HR_G(z)$. Then, using a corollary
to Weyl's inequality (cf. \cite[Chapter~II, \S~3]{gokr1}) we get
\begin{eqnarray*}
|D(re^{i\theta})| 
&\le& \prod_k \left( 1+r^2|\mu_k(HR_G)|^2 \,\right) \\
&\le& \prod_k \left( 1+r^2s_k^2(HR_G)\,\right) \\
&\le& \prod_k \left( 1+r^2||R_G(re^{i\theta})||^2\,s_k^2(H)\,\right) \\
&\le& \prod_k \left(1+\frac{r^2}{|\sin\theta|^2}\,s_k^2(H)\,\right)\,,
\qquad \qquad \qquad \qquad \qquad \quad (3.9)
\end{eqnarray*}
since $s_k(AB)\le ||A||s_k(B)$,
and $||R_G(re^{i\theta})||\le |\sin\theta|^{-1}$.

Let $\nu(t;H)$ be the counting function of the sequence 
$\left\{s_k^{-1}(H)\right\}$:
$$
\nu(t;H)={\rm card }\, \left\{k:\, \frac{1}{s_k(H)}\le t \right\}\,,
$$
and let $\eta=r/|\sin\theta|$. Then, continuing estimate (3.9), we
obtain
$$
\log |D(re^{i\theta})| \le \int_0^\infty \log\left(1+\frac{\eta^2}{t^2}
\right) \, d\nu(t;H) = 2\eta^2 \int_0^\infty
\frac{\nu(t;H)}{t(t^2+\eta^2)}\,dt\,. 
\eqno (3.10)
$$
We set
$$
T(\eta) = \eta^2 \int_0^\infty \frac{\nu(t;H)}{t(t^2+\eta^2)}\,dt\,.
$$
This is an increasing function of $\eta$, such that
$$
T(+0)=T'(+0)=0\,,
\eqno (3.11)
$$
and
\begin{eqnarray*}
\int_0^\infty \frac{T(\eta)}{\eta^{p+1}}\,d\eta &=&
\int_0^\infty \frac{\nu(t;H)}{t}\,dt\, \int_0^\infty
\frac{d\eta}{\eta^{p-1}(t^2+\eta^2)} \\ \\
&=& \left( \int_0^\infty \frac{ds}{s^{p-1}(s^2+1)}\right) \, 
\int_0^\infty \frac{\nu(t;H)}{t^{p+1}}\,dt \\  
&=& \frac{\pi}{2p\sin\frac{\pi p}{2}}\, ||H||_p^p\,.
\qquad \qquad \qquad \qquad \qquad \qquad \quad  (3.12)
\end{eqnarray*}

Now, using factorization (3.8) and the upper bound (3.10), we obtain the
lower bound for the function $C_G(z)$:
$$
\log^-|C_G(re^{i\theta})| \le \frac{\log^-|C_A(re^{i\theta})| +
\log^-|C_{A^*}(re^{i\theta})|}{2}
+ T\left(\frac{r}{|\sin\theta|}\right).
\eqno (3.13)
$$
Inequality (3.13) provides us with a lower bound for the
function $\log|C_G(z)|$.
Later, we perform it into an upper bound for that function, 
but first, we must prepare estimates for the RHS of
(3.13). For the second term, we already have estimate (3.12). In order to
estimate the first sum in the RHS of (3.13), we use Theorem~5. 

According to Theorem~6.1 from \cite[Chapter~II]{gokr1},
$$
\gamma_p(C_A) = \sum_k |{\rm Im}\mu_k(A)|^p \le \sum_k |\mu_k(H)|^p 
= ||H||_p^p\,.
$$
Therefore, our Theorem~5 is applicable to the function $C_A(z)$, and
$$
\int_0^\infty \frac{m(r,C_A)}{r^{p+1}}\,dr \le 
\int_0^\infty \frac{\log M(r,C_A)}{r^{p+1}}\,dr \le
K_p\, \left(||H||_p^p+\alpha_p(C_A^{-1}) \right)\,, 
$$
where 
$$
m(r,f)=\frac{1}{2\pi} \int_0^{2\pi} \log^+|f(re^{i\theta})|\, d\theta
$$ 
is the Nevanlinna proximity function. Since $C_A(0)=1$, making use of
Jensen's formula, we obtain $m(r,C_A^{-1})\le m(r,C_A)$, and
therefore 
$$
\int_0^\infty \frac{m(r,C_A^{-1})}{r^{p+1}}\,dr \le 
\int_0^\infty \frac{m(r,C_A)}{r^{p+1}}\,dr \le 
K_p\, \left(||H||_p^p+\alpha_p(C_A^{-1}) \right)\,. 
\eqno (3.14)
$$
Since zeros of $C_{A^*}$ are conjugates to those of $C_A$, 
$m(r,C_{A^*}^{-1})=m(r,C_A^{-1})$, 
so we have the same upper bound for $m(r,C_{A^*}^{-1})$. 

\begin{procl}{Claim} Estimate (3.13) together with bounds (3.10), (3.12)
and
(3.14) yield the upper bound
$$
\int_0^\infty \frac{\log|C_G(iy)|}{y^{p+1}}\, dy \le 
K_p\, \left(||H||_p^p+\alpha_p(C_A^{-1}) \right)\,. 
\eqno (3.15)
$$
\end{procl}

First, assuming (3.15), we complete the proof of the theorem, and then we
shall prove the claim. 

Since $C_G$ is a canonical product of genus one, we have for ${\Im}z\ne
0$
$$
\log|C_G(z)| = {\rm Re}\, \int_{-\infty}^\infty
\left[\log\left(1-\frac{z}{t}\right) +\frac{z}{t}\right]\,dh_G(t)\,,
$$
where $h_G(t)$ equals the number of zeros of $C_G$ on $(0,t]$, if $t>0$,
and equals the number of zeros of $C_G$ on $[t,0)$ with the minus sign, if
$t<0$. Integrating by parts in the RHS, we obtain 
$$
\log|C_G(z)| = {\rm Re}\,\left[ z^2 \int_{-\infty}^\infty
\frac{h_G(t)}{t^2(z-t)}\,dt \right]\,,
$$
whence 
$$
\log|C_G(iy)| = {\rm Re}\,\left[ -y^2 \int_{-\infty}^\infty
\frac{h_G(t)}{t^2(iy-t)}\,dt \right] = 
y^2 \int_{-\infty}^\infty
\frac{h_G(t)}{t(t^2+y^2)}\,dt\,.
$$

Now, we integrate the last relation by $y$ and apply Fubini's theorem (the
integrand is non-negative):
\begin{eqnarray*}
\int_0^\infty \frac{\log|C_G(iy)|}{y^{p+1}}\,dy
&=&
\int_{-\infty}^\infty \frac{h_G(t)}{t}\,dt\,\int_0^\infty
\frac{dy}{y^{p-1}(t^2+y^2)} \\ \\  
&=& \frac{\pi}{2\sin\frac{\pi p}{2}}\, \int_{-\infty}^\infty
\frac{h_G(t)}{t|t|^p}\,dt \,.
\qquad \qquad \qquad \qquad \quad  (3.16)
\end{eqnarray*}
Combining (3.15) and (3.16) we obtain 
$$
\int_{-\infty}^\infty \frac{h_G(t)}{t|t|^p}\,dt = 
\frac{2\sin\frac{\pi p}{2}}{\pi}  
\int_0^\infty \frac{\log|C_G(iy)|}{y^{p+1}}\,dy 
\le K_p\,\left(||H||_p^p + \alpha_p(C_A^{-1}) \right)\,,
$$
whence
$$
||G||_p^p = \sum_k |\mu_k (G)|^p = \int_{-\infty}^\infty
\frac{dh_G(t)}{|t|^p} 
= p\,\int_{-\infty}^\infty \frac{h_G(t)}{t|t|^p}\,dt  
\le K_p\,\left(||H||_p^p + \alpha_p(C_A^{-1}) \right)\,,
$$
or $||G||_p \le K_p\,\left(||H||_p+\alpha_p(C_A^{-1}) \right)$.
Since $||A||_p\le ||G||_p +||H||_p$, we are done. 

\medskip
It remains to prove the claim. In the proof we use the following
\begin{procl}{Lemma}
Let $u(z)$ be a function harmonic in the upper half-plane, continuous up
to the real axis, and such that
$$
|u(z)|=o(|z|)\,, \qquad z\to 0\,.
\eqno (3.17)
$$
Then
$$
u(re^{i\theta}) \le 
\frac{K}{\sin\theta}\, 
\left\{\frac{1}{\pi} \int_0^{\pi} u^-(2re^{i\varphi})\sin\varphi\,
d\varphi + \frac{r}{2\pi} \int_{-2r}^{2r} \frac{u^-(t)}{t^2}\,dt
\right\}\,.
\eqno (3.18)
$$
\end{procl}

Since the proof of the lemma is a standard
combination of R.~Nevanlinna's representation of functions harmonic in the
upper semi-disk and of Carleman's formula (cf. \cite[Chapter~I]{matsaev2}
and \cite[Lecture~23]{levin}), we give it in the appendix. 

\medskip\par\noindent{\it Proof of the Claim:} We fix $\delta$,
$1<\delta<p$, set $\beta=(\pi-\pi/\delta)/2$, and define the function
$$
u(z)=\log|C_G(z^{1/\delta}e^{i\beta})|\,.
$$
This function is harmonic within the angle 
$$
(1-\delta)\frac{\pi}{2}<\arg z<(1+\delta)\frac{\pi}{2}\,,
$$
and, since $\log|C_G(z)|=O(|z|^2)$, $z\to 0$, it  satisfies condition
(3.17). We take an arbitrary $\psi$ such that $|\psi|\le
(\delta-1)\pi/4$
and apply estimate (3.18) to the function $u_\psi(z)=u(ze^{i\psi})$ 
choosing $\theta=\arg z=\pi/2-\psi$. This way, 
we obtain an estimate for $\log|G(iy)|=u(iy^\delta)$:
$$
\log|C_G(iy)| \le K\, \left\{
\frac{1}{\pi}\int_0^\pi u^-\left(2y^\delta
e^{i(\varphi+\psi)}\right)\sin\varphi \,d\varphi 
+ y^\delta \int_{-2y^\delta}^{2y^\delta}
\frac{u^-(te^{i\psi})}{t^2}\, dt
\right\}\,.
$$
Integrating by $\psi$ and increasing, if needed, the intervals of
integrations, we obtain
$$
\log|C_G(iy)| \le 
K_p\, \left\{\frac{1}{\pi}
\int_{I}
u^-\left(2y^\delta
e^{i\psi}\right) \,d\psi 
+ y^\delta \int_0^{2y^\delta} \frac{dt}{t^2}
\int_{I}
u^-(te^{i\psi})\,d\psi \right\}\,,
$$
where $I=[-(\delta-1)\pi/4,\pi+(\delta-1)\pi/4]$. 
Since
$$
\int_I 
u^-(te^{i\psi}) \,d\psi 
= \int_J \log^-|C_G(t^{1/\delta}e^{i\eta})|\,
d\eta\,,  
\eqno (3.19)
$$
where $J=[\beta_p,\pi-\beta_p]$ with $\beta_p>0$, we obtain
$$
\log|C_G(iy)| \le 
K_p\left\{\int_J
\log^-|C_G(2^{1/\delta}ye^{i\eta} )|d\eta 
+ y^\delta \int_0^{2y^\delta} \frac{dt}{t^2}
\int_J \log^-|C_G(t^{1/\delta}e^{i\eta})| d\eta \right\}.
$$
Plugging here estimate (3.13), and using monotonicity of
$T(\eta)$,
we obtain
\begin{eqnarray*}
\log|C_G(iy)| &\le& 
K_p\, \left\{
m\left(2^{1/\delta}y,C_A^{-1}\right) +
T\left(\frac{2^{1/\delta}y}{\sin\beta_p} \right) 
\right. \\ 
&+& \left. y^\delta \int_0^{2y^\delta} \frac{dt}{t^2}
\left[
m\left(t^{1/\delta},C_A^{-1}\right) +
T\left(\frac{t^{1/\delta}}{\sin\beta_p} \right) 
\right]
\right\}\,.
\qquad \qquad \quad (3.20)
\end{eqnarray*}
Integrating (3.20) with respect to $y$, we get
$$
\int_0^\infty \frac{\log|C_G(iy)|}{y^{p+1}}\,dy 
\le K_p\, \left\{
\int_0^\infty \frac{m(y,C_A^{-1})}{y^{p+1}}\,dy +
\int_0^\infty \frac{T(y)}{y^{p+1}}\,dy
\right. 
$$
$$
+ \left. \int_0^\infty \frac{dy}{y^{1+p-\delta}}
\int_0^{2y^\delta} \frac{dt}{t^2}
\left[
m\left(t^{1/\delta},C_A^{-1}\right) +
T\left(\frac{t^{1/\delta}}{\sin\beta_p} \right) 
\right]
\right\}\,.
\eqno (3.21)
$$
Changing the order of integration, we obtain
$$
\int_0^\infty \frac{dy}{y^{1+p-\delta}} \, \int_0^{2y^\delta}
\frac{m(t^{1/\delta},C_A^{-1})}{t^2}\, dt
$$
$$
= \int_0^\infty \frac{m(s,C_A^{-1})}{s^{\delta+1}}\,ds \,
\int_{s2^{-1/\delta}}^\infty \frac{dy}{y^{1+p-\delta}} 
\le K_p \int_0^\infty \frac{m(s,C_A^{-1})}{s^{p+1}}\,ds \,,
\eqno (3.22)
$$
and similarly
$$
\int_0^\infty \frac{dy}{y^{1+p-\delta}}\, \int_0^{2y^\delta}
\frac{dt}{t^2}\, T\left( \frac{t^{1/\delta}}{\sin\beta_p}\right)
\le K_p\,\int_0^\infty \frac{T(s)}{s^{p+1}}\,ds\,. 
\eqno (3.23)
$$

At last, plugging estimates (3.22) and (3.23) into (3.21), and using
inequalities (3.11) and (3.13), we obtain (3.14), and complete the proof
of the claim, and therefore, of Theorem~4. $\Box$. 
 
\medskip
Observe that the same method based on the factorization (3.8) proves the
following result which is probably known to the specialists:
\begin{procl}{Theorem 8} 
Let $A$ be a compact operator, and let $\mu (A) = \{\mu_k(A)\}$ be its
sequence of eigenvalues counted with multiplicities. Them, for $1<p\le
2$, 
$$
||A||_p\le K_p\,\left( ||H||_p+||\mu (A)||_{l^p} \right)\,.
$$
\end{procl}

\subsection{Appendix. Proof of the Lemma}

According to the Nevanlinna-Green formula (cf.
\cite[Section~24.3]{levin}), for any harmonic function in
the semi-disk $\{{\rm Im}z\ge 0, \, |z|<R\}$ which is continuous up to the
boundary,
$$
u(z)=
\frac{R^2-|z|^2}{2\pi} \,
\int_0^\pi \left(
\frac{1}{|Re^{i\varphi}-z|^2}-\frac{1}{|Re^{i\varphi}-{\bar z}|^2} 
\right) u(re^{i\varphi})\, d\varphi
$$
$$
+\frac{{\rm Im}z}{\pi}\, \int_{-R}^R
\left( \frac{1}{|t-z|^2} - \frac{R^2}{|R^2-tz|^2}
\right) u(t)\,dt
$$
$$
\le 
\frac{R^2-|z|^2}{2\pi} \,
\int_0^\pi \left(
\frac{1}{|Re^{i\varphi}-z|^2}-\frac{1}{|Re^{i\varphi}-{\bar z}|^2} 
\right) u^+(re^{i\varphi})\,d\varphi
$$
$$
+\frac{{\rm Im}z}{\pi}\, \int_{-R}^R
\left( \frac{1}{|t-z|^2} - \frac{R^2}{|R^2-tz|^2}
\right) u^+(t)\,dt
\eqno (A1)
$$
since the both kernels are non-negative. Now, we estimate the kernels for
$z=re^{i\theta}$ and $R=2r$. For the first kernel $K_1(Re^{i\varphi},z)$
we have 
$$
K_1(Re^{i\varphi},z) \le 
\frac{R^2-r^2}{2\pi}\, \frac{4Rr\sin\theta \sin\varphi}{(R-r)^4}
\le \frac{12\sin\varphi}{\pi}\,;
\eqno (A2)
$$
and for the second kernel $K_2(t,z)$ we obtain
$$
K_2(t,z) \le \frac{r\sin\theta}{\pi}\,  
\frac{R^2(R^2-t^2)}{|t-z|^2|R^2-tz|^2} 
$$
$$
\le \frac{r\sin\theta}{\pi}\, \frac{4(R^2-t^2)}{t^2 R^2 \sin^2\theta} 
\,\left(1-\frac{|t|}{2R} \right)^{-2} 
\le \frac{48r}{\pi\sin\theta} \left(\frac{1}{t^2}-\frac{1}{R^2}
\right)\,.
\eqno (A3)
$$
Plugging estimates (A2) and (A3) into (A1), we obtain
$$
u(z)\le \frac{12}{\pi}
\int_0^\pi u^+(Re^{i\varphi})\sin\varphi\,d\varphi
+ \frac{48r}{\pi\sin\theta}
\int_{-R}^R \left(\frac{1}{t^2}-\frac{1}{R^2} \right)u^+(t)\,dt\,.
\eqno (A4)
$$

Now, we use Carleman's formula (cf. \cite[Lecture~24]{levin}). Due to
condition (3.17), it has no error term, and therefore
$$
\frac{1}{\pi R}
\int_0^\pi u^+(Re^{i\varphi})\sin\varphi\,d\varphi
+ \frac{1}{2\pi}
\int_{-R}^R \left(\frac{1}{t^2}-\frac{1}{R^2} \right)u^+(t)\,dt
$$
$$
=
\frac{1}{\pi R}
\int_0^\pi u^-(Re^{i\varphi})\sin\varphi\,d\varphi
+ \frac{1}{2\pi}
\int_{-R}^R \left(\frac{1}{t^2}-\frac{1}{R^2} \right)u^-(t)\,dt\,.
\eqno (A5)
$$
Combining (A4) and (A5), we obtain (3.18) and complete the proof. $\Box$

\bigskip\par\noindent{\em School of Mathematical Sciences, Tel-Aviv
University, \newline
\noindent Ramat-Aviv, 69978, Israel

\medskip\par\noindent matsaev@math.tau.ac.il \newline
\noindent sodin@math.tau.ac.il}

\end{document}